\documentclass[12pt,a4paper]{article} 

\usepackage{amsmath,amssymb,graphicx,color}
\usepackage[english]{babel}
\usepackage{textcomp}
\newtheorem{theoreme}{Theorem}[section]
\newtheorem{lemma}[theoreme]{Lemma}
\newtheorem{define}[theoreme]{Definition}
\newtheorem{propal}[theoreme]{Proposition}

\newenvironment{proof2}
    {\begin{list}{{\bf Proof}}
                 {\setlength{\labelwidth}{0cm}       
                  \setlength{\leftmargin}{0em}}    
                 \item}
    {$\blacksquare$\end{list}}

\def \disp {\displaystyle}

\def \Ch {{\text{\textmarried}^M}}
\def \N {{\mathbb N}} 
\def \R {{\mathbb R}}

\def \Rd {{\mathbb R^3}}
\def \C {{\cal C}}
\def \demi {{\frac {1}{2}}}

\def \eps {{\varepsilon}}
\def \s {{\sigma}}
\def \k {{\kappa}}

\def \tNP {{\tilde{\cal N}}}
\def \ttNP {{\tilde{\tilde{\cal N}}}}
\def \oscil {{w}}
\def \M {{\cal M}}
\def \A {{\cal A}}
\def \D {{\cal D}}
\def \L {{\Lambda}}
\def \E {{{\cal E}}}
\def \En {{{\cal E}_{n}^{\ell,\by}}}

\def \un {{\ensuremath{\textrm{\textnormal{1}}\hspace{-0.7ex}
                       \textrm{\textnormal{I}}}}}

\def \bx {{\bf x}}  \def \nx {{\breve{x}}} 
\def \bX {{\bf X}}  \def \nX {{\breve{X}}} 
                    \def \nB {{\breve{B}}} 
\def \bn {{\bf n}}   
\def \bl {{\bf l}}
\def \by {{\bf y}}  \def \ny {{\breve{y}}} 
\def \bv {{\bf v}}  \def \nv {{\breve{v}}}  
\def \bw {{\bf w}}  \def \nw {{\breve{w}}} 
\def \bz {{\bf z}}  \def \nz {{\breve{z}}}  
\def \bs {{{\boldsymbol\sigma}}}  
                    \def \nW {{\breve{W}}}

\def \ZLy {{Z^{\L,\by}}}
\def \rm {{r_-}}
\def \rp {{r_+}}

\author{
Myriam \textsc{Fradon}\\
{\footnotesize UMR CNRS 8524, UFR de Math\'ematiques}\\
{\footnotesize Universit\'e de Lille 1}\\
{\footnotesize 59655 Villeneuve d'Ascq Cedex, France}\\
{\footnotesize e-mail~:~Myriam.Fradon@univ-lille1.fr}\\
{\footnotesize tel~:~+33(0)32043 6694, fax~:~+33(0)32043 4302}\\
\and 
Sylvie \textsc{R\oe lly}\\
{\footnotesize Institut f\"ur Mathematik der Universit\"at Potsdam}\\ 
{\footnotesize Am Neuen Palais, 10}\\
{\footnotesize 14469 Potsdam, Germany}\\
{\footnotesize e-mail~:~roelly@math.uni-potsdam.de}\\
{\footnotesize tel~:~+49 (0)331 9771478, fax~:~+49 (0)331 9771001}\\
{\footnotesize On leave of absence Centre de
  Math\'ematiques Appliqu\'ees,} \\
{\footnotesize UMR C.N.R.S. 7641, \'Ecole
  Polytechnique,}\\
{\footnotesize 91128 Palaiseau Cedex, France.}
  }
\date{}

\title{Infinitely many Brownian globules with Brownian radii}

\begin{document}

\maketitle

\vspace{3cm}

\centerline{\textbf{Abstract}}
We consider an infinite system of non overlapping  globules undergoing Brownian motions in $\R^3$.
The term {\it globules} means that the objects we are dealing with are spherical, but with a radius which is random and time-dependent.
The dynamics is modelized by an infinite-dimensional Stochastic Differential Equation
with local time.
Existence and uniqueness of a strong solution is proven for such an equation 
with  fixed deterministic initial condition. We also find a class of reversible measures.\\

\noindent
AMS Classifications: 60H10, 60J55, 60J60.\\
KEY-WORDS: Infinite-dimensional Stochastic Differential Equation, hard core potential, oblique reflection, reversible measure, local time.\\

\newpage

\section{Introduction}

The aim of this paper is to construct a random dynamics performed by an infinite system of globules, where a globule is a sphere in $\R^3$ with variable  radius. The centers of the globules undergo independent Brownian motions, while their radii perform Brownian oscillations between a minimum and a maximum value. Since the scale of these oscillations can be different than those of the centers, we introduce a coefficient  $\s$ which reflects the elasticity of the surface of each globule. The globules  can not overlap and when the distance between two globules becomes 0, they repel each other immediately; that means they interact through a hard core potential.

A reversible system of infinitely many Brownian hard spheres (called balls) was first introduced and analyzed by H. Tanemura \cite{TanemuraEDS}. 
Then, some natural generalizations were studied for different types of additional smooth interactions between the balls:  for a gradient type interaction with finite range in \cite{FR1},\cite{FR2} and  for an interaction with infinite range in \cite{FRT}. 
The specificity one has to deal with in a hard core situation -- hard balls can not overlap -- comes from the additional infinite-dimensional local time term in the SDE. 
Notice that in all these works the spheres have a fixed positive radius. 

The originality and new difficulty of the present model -- which can find relevant applications in cell dynamics like molecular motors -- lies in the random oscillations of the radius of each sphere. 
We propose here a pathwise approach for the construction of this infinite-dimensional dynamics, by building a sequence of finite-dimensional approximating processes. 
But already for finitely many globules, the existence of such dynamics is not a simple question.
 Indeed, one of the authors constructed recently in \cite{FradonBrownianGlobules} a finite system of mutually repelling 
Brownian globules. Nevertheless, we need here a non trivial generalization of these results: since the scale $\s$ of the radii oscillations is different than the scale of the center oscillations, the direction of the reflection after a collision between two globules is no more normal as in \cite{FradonBrownianGlobules}. It is an oblique reflection of Brownian motions on a complex non smooth domain, whose existence problem we solve in Proposition \ref{Prop:ExistApprox}.   

In Section \ref{Sect_the_model} we present the model and its dynamics described by the stochastic differential equation ($\E$) 
and we state the results. In Section \ref{Sec_approx}, we show the convergence of the approximations 
 and analyze the limit process. Last, we remark that some kind of hard core Poisson measure is reversible for this dynamics.

\section{The infinite model of mutually repelling globules with Brownian radii}\label{Sect_the_model}

A globule is a sphere in $\R^d$ with a variable radius. For $d=2$ the globules modelize for example the motion of discs on a flat surface or balls floating on a liquid. In this paper, we fix $d=3$ which corresponds to the natural physical case of bubbles in the Euclidean space. Our techniques and results obviously extend to any dimension $d$ larger than $1$.\\

A globule is characterized by a pair $(x,\nx)\in  \R^3 \times \R$. ~$x$ is the position of the center of the globule  and $\nx$ is its radius. \\
We are dealing here with infinitely many indistinguishable globules, thus the state space of the system is included in $\M$, the set of point measures on $\R^3 \times \R$.
A  configuration of globules is a locally finite point measure $\bx=\sum_{i \in J} \delta_{(x_i,\nx_i)}$ on $\R^3 \times \R$, where $(x_i,\nx_i)$ characterizes the $i$-th globule and $J \subset \N$.
For simplicity sake, we will identify any such point measure $\bx$ with its support $\{ (x_i,\nx_i), i\in J\} \subset \Rd \times \R$.

The globules we deal with in this paper can not overlap and their radii are bounded from below (resp. from above) by a constant $\rm>0$ (resp. $\rp >\rm$). 
Thus the exact configuration space $\A_g$ of all {\bf allowed configurations of globules} is the following set:
$$
\A_g=\left\{ \begin{array}{l}
              \bx =\{(x_i,\nx_i), i \in J\} \text{ for some } J\subset \N,                                                \\
              \text{ with } x_i\in \Rd, \nx_i\in [\rm;\rp] \textrm{ and } |x_i-x_j|\ge \nx_i+\nx_j \text{ for } i\neq j
             \end{array} \right\}
$$
Let us notice that this model can not be reduced to a hard core model in $\R^3 \times \R$. Indeed, a hard core condition between globules  would mean that there exists $\rho$ such that $\forall i\neq j, \quad |(x_i,\nx_i)-(x_j,\nx_j)|\ge \rho$, which is equivalent to
$|x_i-x_j|^2 + |\nx_i-\nx_j|^2 \ge \rho^2$. This last inequality is clearly not comparable with the condition $|x_i-x_j|\ge \nx_i+\nx_j $.\\

We will use the notations:
\begin{itemize}
\item
$ B(x,\rho)$ is the closed ball centered in $x \in \Rd$ with radius $\rho$ and by extension,
for any subset $A$ in $\Rd $, we define the $\rho$-neighborhood of $A$ by 
$$
B(A,\rho) := \{ y \in \Rd \textrm{ such that } d(y,A) \le \rho \}.
$$
where $d(y,A)$ denotes the Euclidean distance between $y$ and $A$.
\item
The symbol $|v|$ denotes the Euclidean norm of the vector $v$. \\
We also denote the volume of a subset $A$ of $\Rd$ by $|A|$.
\item
For $A \times I$ a Borel subset of $\R^3 \times \R $, $ N_{A \times I} $ is the counting variable on $ \M $~:
$$
 N_{A \times I}(\bx) = \sharp \{ i \in \N : x_i \in A \textrm{ and } \nx_i \in I \}.
$$
\item
For $ \L $ a Borel subset of $\R^3 \times \R $, $ {\cal B}_\L $ is the $ \sigma$-algebra on $ \M $ generated by the sets $ \{ N_{\L'} = n \} $, $ n \in \N $, $ \L' \subset \L $, $ \L' $ bounded.
\item
We write $ \bx_{\L} =\bx \cap \L $ for the restriction of the configuration $\bx$ to $ \L \subset \R^3 \times \R$,\\
and $ \bx \by $ for the concatenation of configurations $\bx$ and $\by$.
\item
$\pi$ (resp. $\pi_\L$) is the Poisson process on $\R^3 \times \R$ (resp. on $\L$) with intensity measure the Lebesgue measure
$d\by$ (resp. $d\by|_\L$).
\end{itemize}

We define the set $\Pi_g$ of {\bf hard globule Poisson processes} via a 
local density function~:

\begin{define}
A Probability measure $ \mu $ on $ \M $ is a hard globule Poisson process 
if and only if, for each compact subset $ \L \subset \R^3 \times \R $,
$$
\mu( d\bx | {\cal B}_{\L^c} )(\by)
= \frac{1}{\ZLy} \un_{\{\bx_\L\by_{\L^c}\in \A_g\}} ~\pi_\L(d\bx) \textrm{ for } \mu \textrm{-a.e. } \by.
$$
where the so-called partition function $\ZLy$ is the renormalizing constant~: 
$$\ZLy
  =  e^{- | \L |}  
  \left( 1 +
  \sum_{n=1}^{+\infty} ~
  \frac{1}{n!} ~ \int_{\L^n}  
                   \un_{\{\bx_\L\by_{\L^c}\in \A_g\}}
                  \,  dx_1 d\nx_1 \cdots dx_n  d\nx_n
  \right) . 
  $$
\end{define}
At least one hard globule Poisson process exists (generalization of the existence results for hard core Gibbs measure in \cite{Dob}).
It is conjectured  that for $ r_+ $ small enough it is unique, while for $ r_- $ large 
enough phase transition occurs~: $ \Pi_g $ should contain several measures (see e.g. \cite{Georgii}).\\

In order to modelize the random motion of repelling globules with oscillating radii, let us consider a probability space $(\Omega,{\cal F},P)$ endowed with a complete filtration  $\{{\cal F}_t\}_{t \ge 0}$ and  two sequences of ${\cal F}_t$-Brownian motions: $ ( W_i(t) , t \ge 0 )_{i \in \N} $ which are independent $\R^3$-valued Brownian motions and  $ ( \nW_i(t) , t \ge 0 )_{i \in \N} $ which are independent $\R$-valued Brownian motions, independent from the $W_i$'s too. \\
We fix a parameter $\s > 0 $ which measures the scale of the radii oscillations, that is the elasticity of the surface of each globule.\\
We consider the following system of stochastic differential equations with (oblique) reflection~:
$$
(\E_g)
\left\{
\begin{array}{l}
 \textrm{For } i\in \N , t \in [0,1],                                     \\
 \disp
  X_i(t) = X_i(0) +  W_i(t)    
           + \sum_{j\in \N} \int_0^t \frac{X_i(s)-X_j(s)}{\nX_i(s)+\nX_j(s)} dL_{ij}(s) \\
 \disp
 \nX_i(t) = \nX_i(0) + \s \nW_i(t) - \s^2 \sum_{j\in \N} L_{ij}(t) - L_{i+}(t) + L_{i-}(t)\\
 \disp
\textrm{ where the local times }L_{ij}, L_{i+}, L_{i-} \textrm{ satisfy }
L_{ij}(t) = \int_0^t \un_{\{|X_i(s)-X_j(s)|=\nX_i(s)+\nX_j(s)\}} ~dL_{ij}(s)~, \\
\disp
L_{i+}(t) = \int_0^t \un_{\{\nX_i(s)=\rp\}} ~dL_{i+}(s) \quad\textrm{ and }\quad L_{i-}(t) = \int_0^t \un_{\{\nX_i(s)=\rm\}} ~dL_{i-}(s) .
\end{array}
\right.
$$
As usual, the collision local times are non-decreasing $\R^+-$valued continuous processes with bounded variations and satisfy $L_{ij} \equiv L_{ji}$ and $L_{ii} \equiv 0  $.
The starting configuration $\bX(0)=\{(X_i(0),\nX_i(0)), i\in \N \}$ is a point in $\A_g$.

A solution of the system $ (\E_g) $ is a family 
$ ( X_i(t),\nX_i(t) , L_{ij}(t), L_{i+}(t), L_{i-}(t), 0\le t \le 1, i,j \in \N) $  of processes satisfying $ (\E_g)$. \\

Let us interpret the different terms of $(\E_g)$~:
\begin{itemize}
\item
when two globules collide ($|X_i(t)-X_j(t)|=\nX_i(t)+\nX_j(t)$), they are deflated ($\nX_i(t) $ decreases by $ dL_{ij}(t)$) and move away from each other ($X_i(t)$ is submitted to the repulsive force $\frac{X_i(t)-X_j(t)}{\nX_i(t)+\nX_j(t)}$);
\item
when the radius of a globule reaches the maximal value ($\nX_i (t) = r_+$), it is deflated ($\nX_i(t)$ decreases by $dL^+_i (t)$);
\item
when the radius of a globule reaches the minimal value  ($\nX_i(t)=\rm$), it is inflated ($\nX_i (t)$ increases by $dL^-_i(t)$).
\end{itemize}

\begin{theoreme}
\label{ThExistenceSolution}
The stochastic equation $ (\E_g)$ admits a unique solution with values in ${\A_g}$ for any deterministic initial configuration which belongs to a full measure subset $\underline{{\A_g}}$ in $\A_g $. 
\end{theoreme}
\begin{propal}
\label{ThReversibiliteSolution}
If the initial distribution $\mu$ is a hard globule Poisson process and if $\mu (\underline{{\A_g}}) =1$, then the solution of $ (\E_g)$ is time-reversible, that is its law is invariant with respect to the time reversal.
\end{propal}

The next section is devoted to the proofs of these results.

\section{The infinite-dimensional process, constructed by approximation}\label{Sec_approx}
\subsection{The approximating processes}\label{SectApproximatingProcesses}

In this whole subsection, $\ell \in \N^*$  and $\by \in \A_g$ are fixed.
Using a classical penalization method with external configuration $\by$, we construct an approximating process which {\bf essentially stays} in the ball $ B(0,\ell)$ .
This is done by introducing in the dynamics $ (\E_g) $ an additional drift which vanishes in a subset of $B(0,\ell)$ and is strongly repulsive outside of $B(0,\ell)$;  we take as drift the gradient of the $C^2_b$-function (twice differentiable function with bounded derivatives) which is defined on $\R^3\times\R$ by:
$$
\psi^{\ell,\by}(x,\nx)= \psi_1(|x|)+\psi_2(\nx)+\sum_{j:|y_j|>\ell} \psi_3\left( \frac{|x-y_j|}{\nx+\ny_j} \right)
$$
with $\psi_1$, $\psi_2$ and $\psi_3$ non-negative $\C^{\infty}-$functions vanishing respectively on $]-\infty,\ell]$, $[\rm,\rp]$ and $[1,+\infty[$, and increasing rapidly on their supports~: \\
$\psi_1(s)=2s$ for $s \ge \ell+e^{-\ell}$ \\
$\psi_2(s)=\ell s$ for $s \ge \rp+e^{-\ell}$ and $\psi_2(s)=\ell(\rp+\rm-s)$  for $s \le \rm-e^{-\ell}$ \\
$\psi_3(s)=\ell$ for $s \le 1-e^{-\ell}$.

The function $\psi^{\ell,\by}$ satisfies
$$
\psi^{\ell,\by}(x,\nx)=0 \quad\Leftrightarrow \quad x \in B(0,\ell) \textrm{ and } (x,\nx)\by_{B(0,\ell)^c} \in \A_g .
$$
By a slight abuse of notation, $\by_{B(0,\ell)^c} $ denotes the restricted configuration $\cup_{\{j: |y_j|>\ell\}} \{(y_j,\ny_j)\} $.\\
Remark that the functions $\psi^{\ell,\by}$  are so repulsive for large $\ell$ that they satisfy
\begin{equation}\label{HypSurPsile}
\quad\quad\quad\quad 
\sup_{\by \in \A_g} \sum_{\ell=1}^{+\infty} 
\int_{\R^3 \times \R} \un_{\psi^{\ell,\by}(x,\nx)>0}~\exp(-\psi^{\ell,\by}(x,\nx))~dx d\nx \quad < \quad + \infty.
\end{equation}
Let us now define the finite-dimensional dynamics~:
$$
(\En) \quad 
\left\{
\begin{array}{l}
  \forall i \in \{1,\ldots,n\}, \quad \forall t \in [0,1],                           \\
  \disp
  X_i(t) = X_i(0) +W_i(t) -\demi \int_0^t \nabla_{x} \psi^{\ell,\by}( X_i(s),\nX_i(s) ) ds 
                          +\sum_{j=1}^n \int_0^t \frac{X_i(s)-X_j(s)}{\nX_i(s)+\nX_j(s)} dL_{ij}(s)\\
 \disp
 \nX_i(t) = \nX_i(0) +\s \nW_i(t) -\frac{\s^2}{2} \int_0^t \nabla_{\nx} \psi^{\ell,\by}( X_i(s),\nX_i(s) ) ds 
                                  -\s^2 \sum_{j=1}^n L_{ij}(t) -L_{i+}(t) +L_{i-}(t)\\
\disp
\textrm{ where the local times satisfy } L_{ij}(t) = \int_0^t \un_{\{|X_i(s)-X_j(s)|=\nX_i(s)+\nX_j(s)\}} ~dL_{ij}(s)~, \\
\disp
L_{i+}(t) = \int_0^t \un_{\{\nX_i(s)=\rp\}}~dL_{i+}(s) \quad\textrm{ and }\quad L_{i-}(t) = \int_0^t \un_{\{\nX_i(s)=\rm\}}~dL_{i-}(s).
\end{array}
\right.
$$
$(\En)$ is a $n$-dimensional reflected stochastic differential equation. \\
If $\s=1$, that is if the radii oscillations and the center oscillations are on the same scale, the above Skorohod equation contains a normal reflection on the boundary of the set of allowed configurations of $n$ globules. The problem of existence and reversibility of this type of dynamics was recently  solved by one of the authors in \cite{FradonBrownianGlobules}. \\
For $\s \not= 1$ a new difficulty occurs. The physically evident reflection on the boundary of the domain of allowed configurations is now oblique since the radii oscillations have another time scale as the center oscillations: it corresponds to a normal reflection but in an anisotropic configuration space, where the radii coordinates are rescaled by $1/\sigma$. The existence of solution for general SDEs with oblique reflections on nonsmooth domains is a hard problem which is solved in the literature only in some particular cases, like for intersection of smooth bounded domains or for polyhedral domains (see e.g. \cite{DupuisIshii} and \cite{Williams}). Since our model is not covered by these works, we present in the following proposition a suitable existence result. 

\begin{propal} \label{Prop:ExistApprox}
Assume that $\Phi$ is an $\R$-valued $\C^2_b$-function defined on $(\R^3\times \R)^n$.
There exists a unique strong solution to the Skorohod problem
$$
(\E^\Phi_n)\left\{
\begin{array}{l}
\forall i \in \{1,\ldots,n\}, \quad \forall t \in [0,1],                           \\
\disp
\disp
X_i(t) = X_i(0) +W_i(t) -\frac{1}{2} \int_0^t \nabla_{x_i}\Phi(\bX(s)) ds +\sum_{j=1}^n \int_0^t \frac{X_i-X_j}{\nX_i+\nX_j}(s) dL_{ij}(s)   \\
\disp
\nX_i(t) = \nX_i(0) +\s\nW_i(t) -\frac{\s^2}{2} \int_0^t \nabla_{\nx_i}\Phi(\bX(s)) ds -\s^2 \sum_{j=1}^n L_{ij}(t) - L_{i+}(t) + L_{i-}(t)\\
\disp
\textrm{ where the local times satisfy } L_{ij}(t) = \int_0^t \un_{\{|X_i(s)-X_j(s)|=\nX_i(s)+\nX_j(s)\}} ~dL_{ij}(s)~, \\
\disp
L_{i+}(t) = \int_0^t \un_{\{\nX_i(s)=\rp\}} ~dL_{i+}(s) \quad\textrm{ and }\quad L_{i-}(t) = \int_0^t \un_{\{\nX_i(s)=\rm\}} ~dL_{i-}(s) .
\end{array}
\right.
$$ 
For any initial condition in  $\A_g$ the solution is an $\A_g$-valued process. \\
Moreover, the solution with initial distribution $e^{-\Phi(\bx)} \un_{\A_g}(\bx) d\bx$ is time-reversible.
\end{propal}
The proof of this proposition is postponed to the end of this section. 
A key idea is the transformation (see (\ref{equ:ProcTransf})) of the initial Skorohod problem into a simpler one, by stretching the radii coordinates by the factor $1/\sigma$ and thus transforming the oblique reflection in a normal one on a modified domain. Let us underline that the oblique reflection we consider is the unique one for which the existence of a reversible dynamics is ensured.\\

Applying Proposition \ref{Prop:ExistApprox} with the potential $\Phi(\bx) =\sum_{i=1}^n \psi^{\ell,\by}(x_i,\nx_i)$, we obtain the existence of a solution to  $(\En)$. \\
When the initial condition is the deterministic configuration $\bx$, this solution is denoted by $\bX^{\ell, \by,n}(\bx,\cdot)$. 
In particular, the $\A_g$-valued finite-dimensional process $\bX^{\ell,\by}$ with initial configuration $\bx=\by_{B(0,\ell)}$ 
and external configuration $\by_{B(0,\ell)^c}$ evolving under the random dynamics $({\cal E}_{n}^{\ell,\by})$ is~:
$$
\bX^{\ell, \by}(\cdot) := \bX^{\ell, \by,n}(\by_{B(0,\ell)},\cdot) 
\quad \textrm{ with } n=\sharp \{ i \in \N : y_i \in B(0,\ell) \} .
$$
The associated local times are denoted by $ L_{i,j}^{\ell, \by}, L_{i+}^{\ell, \by}, L_{i-}^{\ell, \by}, i,j \in \N$.\\

If the initial condition of the system $(\En)$ is random with distribution given by the finite measure:
$$
\nu_n^{\ell,\by}(d\bx) := 
\exp( - \sum_{i=1}^n \psi^{\ell,\by}(x_i,\nx_i)) ~ \un_{\A_g} (\bx) ~dx_1 d\nx_1 \cdots dx_n  d\nx_n.
$$
then the solution of $(\En)$ is reversible. Its law is denoted by $Q_n^{\ell,\by}$.\\
Consider now the following Poisson mixture in $n$ of the $Q_n^{\ell,\by}$'s:
$$
Q^{\ell,\by} = \frac{e^{-|B(0,\ell)|}}{Z^{\ell,\by}} \sum_{n=0}^{+\infty} \frac{1}{n!} ~ Q_n^{\ell,\by} ,
$$
where $Z^{\ell,\by}$ is the renormalizing constant which ensures that $Q^{\ell,\by}$ is a Probability measure.
As a mixture of $ \A_g$-supported time-reversible measures, $Q^{\ell,\by}$ is time-reversible with support included in $ \A_g $. 
Its projection at time 0 is the Probability measure
$$
\mu^{\ell,\by}(d\bx) := 
\frac{e^{-|B(0,\ell)|}}{Z^{\ell,\by}} \sum_{n=0}^{+\infty} \frac{1}{n!} ~\nu_n^{\ell,\by}(d\bx) ,
$$
which represents the law of a Poissonian number of globules essentially concentrated in $B(0,\ell)$.\\

We will construct the infinite-dimensional globule process as limit in $\ell$ of $\bX^{\ell, \by}$; unfortunately, $\bX^{\ell, \by}$ is not time-reversible. This is why we had to introduce $Q^{\ell,\by}$, whose reversibility plays a crucial role in the study of the set of nice paths defined in the next section. 
Moreover, we will prove that the law of $\bX^{\ell, \by}$ and $Q^{\ell,\by}$ are asymptotically close.

To complete this section, let us prove Proposition \ref{Prop:ExistApprox}.
\begin{proof2}
We first introduce an anisotropic linear transformation $\bs^{-1}$ on the space of globule configurations by 
$$
\bx^\s := \bs^{-1}\bx \qquad \Leftrightarrow \qquad \forall i,\quad  x^\s_i=x_i \textrm{  and } \nx^\s_i=\frac{1}{\s} \nx_i .
$$ 
We also transform the process, the potential and the local times as follows~:
\begin{equation} \label{equ:ProcTransf}
\bX^\s = \bs^{-1}\bX, \quad \Phi^\s(\cdot)=\Phi(\bs\cdot),~~ L^\s_{ij} =\sqrt{2+2\s^2} L_{ij},~~ L^\s_{i+}= \frac{1}{\s} L_{i+},~~ L^\s_{i-}=\frac{1}{\s} L_{i-}.
\end{equation}
Moreover, the set of allowed configurations becomes
$$
\A_g^{\s}:=\{ \bx: \bs\bx \in \A_g \} .
$$
$\bX$ and its associated local times are solution of the system  $(\E^\Phi_n)$ if and only if $\bX^\s$ and its associated transformed local times are solution of the following system:
$$
(\E^{\Phi,\s}_n)\left\{
\begin{array}{l}
 \disp
  X^\s_i(t) = X^\s_i(0) +W_i(t) -\frac{1}{2} \int_0^t \nabla_{x_i}\Phi^\s(\bX^\s(s)) ds
              + \sum_{j=1}^n \int_0^t \frac{1}{\sqrt{2+2\s^2}}\frac{X^\s_i-X^\s_j}{\s(\nX^\s_i+\nX^\s_j)}(s) dL^\s_{ij}(s) \\
 \disp
 \nX^\s_i(t) = \nX^\s_i(0) +\nW_i(t) -\frac{1}{2} \int_0^t \nabla_{\nx_i}\Phi^\s(\bX^\s(s)) ds 
               -\frac{\s}{\sqrt{2+2\s^2}} \sum_{j=1}^n L^\s_{ij}(t) - L^\s_{i+}(t) + L^\s_{i-}(t) \\
\disp
\text{where the local times satisfy } 
L^\s_{ij}(t) = L^\s_{ji}(t) = \int_0^t \un_{\{|X^\s_i(s)-X^\s_j(s)|=\s(\nX^\s_i(s)+\nX^\s_j(s))\}} ~dL^\s_{ij}(s)~, \\
\disp
L^\s_{i+}(t) = \int_0^t \un_{\{\nX^\s_i(s)=\frac{\rp}{\s}\}} ~dL^\s_{i+}(s) \quad\text{ and }\quad 
L^\s_{i-}(t) = \int_0^t \un_{\{\nX^\s_i(s)=\frac{\rm}{\s}\}} ~dL^\s_{i-}(s) .
\end{array}
\right.
$$ 
Furthermore, the reversibility of $\bX$ is equivalent to the reversibility of $\bX^\s $; with other words, the solution of $(\E^\Phi_n)$ with initial distribution $\frac{1}{Z}e^{-\Phi(\bx)} \un_{\A_g}(\bx) d\bx$ is reversible if and only if the solution of $(\E^{\Phi,\s}_n)$ with initial distribution $\frac{1}{Z^\s}e^{-\Phi^\s(\by)} \un_{\A_g^\s}(\by) d\by$ is reversible.

The new system of globules  $(\E^{\Phi,\s}_n)$ has now the form of a Skorohod problem with normal reflection. Thus it has a unique solution under the assumptions of Theorem 3.3 (and Corollary 3.6) in \cite{FradonBrownianGlobules}, that is if the domain $\A_g^\s$ on which the equation is reflected satisfies the geometrical regularity properties listed in \cite{FradonBrownianGlobules} Proposition 3.4. The rest of the proof consists in showing these four properties (see \cite{FradonBrownianGlobules} for the relevant definitions).

The domain $\A_g^\s$ does not have a smooth boundary but it is the intersection of smooth domains in the following way: 
$$
\A_g^{\s}=(\bigcap_{1 \le i<j \le n} \D_{ij}) \cap (\bigcap_{1 \le i \le n} \D_{i+}) 
                                              \cap (\bigcap_{1 \le i  \le n} \D_{i-})
$$ 
where $\disp \D_{ij}=\left\{ \bx,~ |x_i-x_j|\ge \s(\nx_i+\nx_j) \right\}$,
$\disp \D_{i+}=\left\{ \bx,~ \nx_i \le \frac{\rp}{\s} \right\}$ and $\disp \D_{i-}=\left\{ \bx,~ \nx_i \ge \frac{\rm}{\s} \right\}$.

{\it (i) At each point $\bx$ of the boundary of the smooth set $\D_{ij}$ (resp. $\D_{i+}$, $\D_{i-}$), there exists a unique unit normal vector $\bn_{ij}(\bx)$ (resp. $\bn_{i+}$, $\bn_{i-}$)}.

Each $\D_{ij}$ is a smooth set with unit inward normal vector at point $\bx\in\partial\D_{ij}$ equal to 
$$
\bn_{ij}(\bx)=\frac{\bw}{\sqrt{2+2\s^2}} \textrm{ where }
w_i=\frac{x_i-x_j}{\s(\nx_i+\nx_j)}=-w_j, \nw_i=\nw_j=-\s \textrm{ and } w_k=\nw_k=0, k\not =i,j. 
$$
Each $\D_{i+}$ (respectively $\D_{i-}$) is a  half-space of $(\R^3 \times \R)^n$ with a constant unit inward normal vector
$\bn_{i+}=\big(0,\cdots,0,-1,0,\cdots,0\big) \, \text{ ($(2i-1)^\text{th}$ coordinate equal to $-1$)}$  (resp. $\bn_{i-}= -\bn_{i+} $).\\
\definecolor{VertNoel}{rgb}{0,0.5,0}
\begin{figure}[h!]
  \centering
  \includegraphics[scale=0.9]{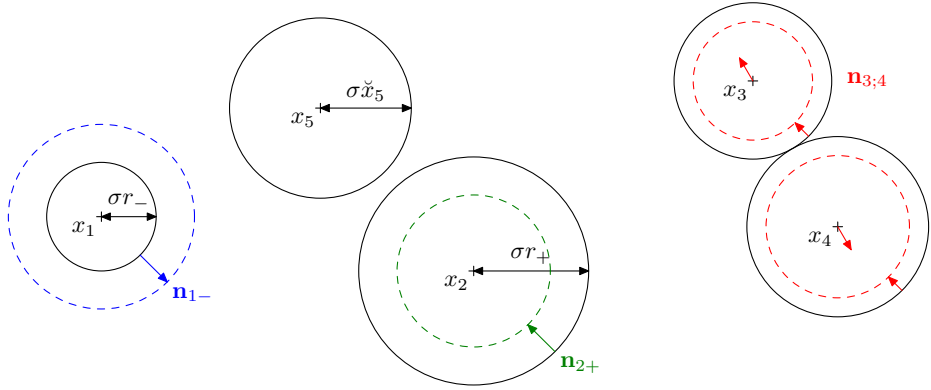}
  \caption{A configuration of 5 globules in $\partial\D_{1-} \cap \partial\D_{2+} \cap \partial\D_{3;4} \cap \A_g^\s$, and the different directions  of impulsion {\color{blue}$\bn_{1-}$}, {\color{VertNoel}$\bn_{2+}$} or {\color{red}$\bn_{3;4}(\bx)$} to go back into the interior of $\A_g^\s$.}
\end{figure}

{\it (ii) Each set $\D_{ij}$ has the Uniform Exterior Sphere property on $\A_g^{\s}$:} \\
$$
\exists \alpha_{ij}>0 , \quad \forall \bx\in \A_g^{\s}\cap\partial\D_{ij} \quad \mathring{B}(\bx-\alpha_{ij}\bn_{ij}(\bx),\alpha_{ij}) \cap \D_{ij} = \emptyset .
$$
\\
Each $\bx\in\A_g^{\s}\cap\partial\D_{ij}$ satisfies $|x_i-x_j|=\s(\nx_i+\nx_j)\ge 2\rm$. \\
For $\bx^{es}=\bx-\rm\sqrt{2+2\s^2}\bn_{ij}(\bx)$ and for any $\bz$~:
$$
\begin{array}{l}
|(x^{es}_i+z_i)-(x^{es}_j+z_j)|-\s(\nx^{es}_i+\nz_i+\nx^{es}_j+\nz_j)                                                                            \\ 
\disp
\quad\le |z_i|+|z_j|-\s(\nz_i+\nz_j)+\left| x_i-x_j-2\rm\frac{x_i-x_j}{\s(\nx_i+\nx_j)} \right|-\s(\nx_i+\nx_j)-2\rm\s^2     \\
\quad\le \sqrt{2+2\s^2}|\bz|-2\rm(1+\s^2) 
\end{array}
$$
This is negative as soon as $\disp |\bz|<\rm\sqrt{2+2\s^2}$. 
Consequently, the property {\it (ii)} holds with $\alpha_{ij} \equiv \rm\sqrt{2+2\s^2}$. See figure 2.
\definecolor{Brun}{rgb}{0.27,0.14,0.04}
\begin{figure}[h!] \label{figure2}
  \centering
  \includegraphics[scale=0.8]{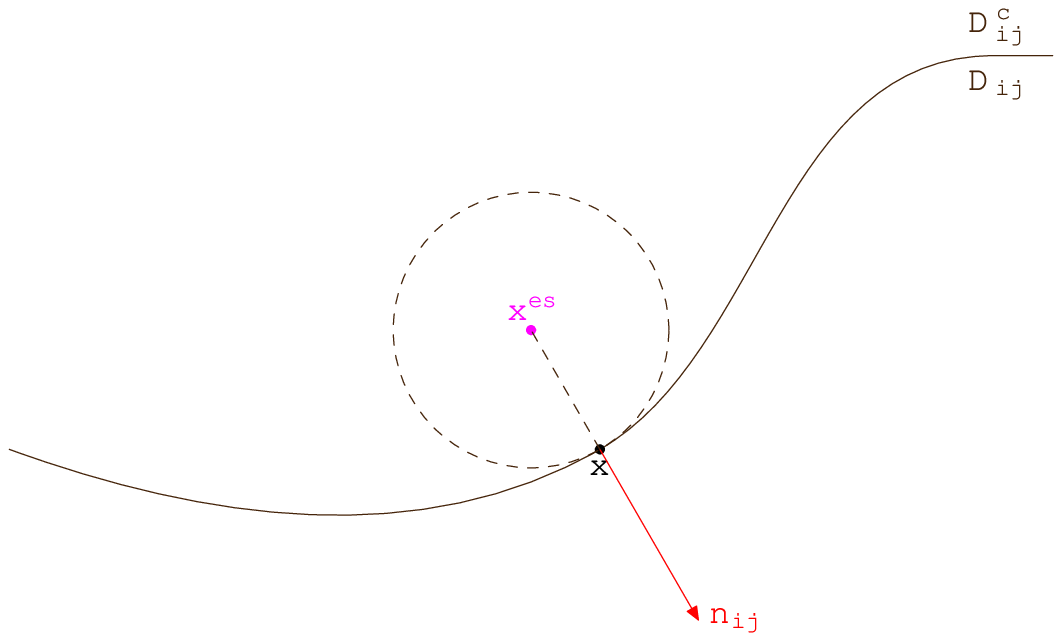}
  \includegraphics[scale=1]{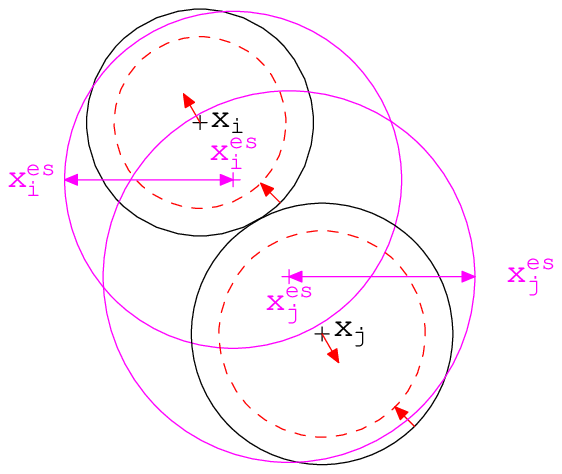}
  \caption{A configuration $\bx\in\partial\D_{ij}$ and the corresponding center {\color{magenta}$\bx^{es}$} of the Uniform Exterior Sphere : {\color{magenta} $\bx^{es}=\bx - \alpha_{ij} \bn_{ij}(\bx) = \bx-\rm\sqrt{2+2\s^2}\bn_{ij}(\bx)$}.  Left, a simplified representation in $\R^2$. Right, a representation as a pair of colliding globules. Remark that $\bx^{es} \not\in \A_g^{\s}$.}
\end{figure}

{\it (iii) Each set $\D_{ij}$ has the Uniform Normal Cone property on $\A_g^{\s}$: \\
$\exists \beta_{ij}\in [0,1[$ and $\delta_{ij}>0$ such that, for each $\bx\in\A_g^{\s}\cap\partial\D_{ij}$, there is a unit vector $\bl^{ij}_\bx$ satisfying
$$ 
\by\in\A_g^{\s}\cap\partial\D_{ij} \cap B(\bx,\delta_{ij}) \Rightarrow  \bn_{ij}(\by).\bl^{ij}_\bx \ge \sqrt{1-\beta_{ij}^2}.
$$}
For $\bx,\by\in\A_g^{\s}\cap\partial\D_{ij}$ one has 
$\disp \bn_{ij}(\by).\bn_{ij}(\bx)=\frac{1}{1+\s^2}\left(\frac{x_i-x_j}{|x_i-x_j|}.\frac{y_i-y_j}{|y_i-y_j|}+\s^2\right)$.
The inequality
$$
\frac{x_i-x_j}{|x_i-x_j|}.\frac{y_i-y_j}{|y_i-y_j|}
\ge \frac{1-\sqrt{2}\frac{|\bx-\by|}{|x_i-x_j|}}{1+\sqrt{2}\frac{|\bx-\by|}{|x_i-x_j|}}
\ge 1-2\sqrt{2}\frac{|\bx-\by|}{|x_i-x_j|}
\ge 1-\frac{\sqrt{2}}{\rm}|\bx-\by|
$$
implies that $\disp \bn_{ij}(\by). \bn_{ij}(\bx) \ge 1-\frac{\sqrt{2}}{\rm(1+\s^2)}|\bx-\by|$. 
Hence, for $\bl^{ij}_\bx=\bn_{ij}(\bx)$, {\it (iii)} is satisfied for each $\beta_{ij}\in [0,1[ $ as soon as 
$\disp \delta_{ij} \le \frac{\rm(1+\s^2)}{\sqrt{2}}(1-\sqrt{1-\beta_{ij}^2})$.

{\it (iv) Compatibility between the boundaries: \\
$\exists \beta_0>\sqrt{2\max_{i,j} \beta_{ij}},
\, \forall\bx\in\partial\A_g^{\s}, \, \exists \, \bv(\bx) \text{ such that } \left\{
\begin{array}{l}
\bx\in\partial\D_{ij} \Rightarrow  \bv(\bx).\bn_{ij}(\bx) \ge \beta_0 |\bv(\bx)|\\ 
\bx\in\partial\D_{i+} \Rightarrow  \bv(\bx).\bn_{i+}(\bx) \ge \beta_0 |\bv(\bx)|\\
\bx\in\partial\D_{i-} \Rightarrow  \bv(\bx).\bn_{i-}(\bx) \ge \beta_0 |\bv(\bx)|
\end{array}
\right.$.}

Let us define the following cluster
$$
C_i(\bx)= \{ i \} \cup \left\{ k_m,~\exists k_1,\ldots,k_m \text{ such that } \bx\in\partial\D_{i k_1} \cap 
                                       \partial\D_{k_1 k_2} \cap\cdots\cap \partial\D_{k_{m-1} k_m}  \right\}
$$
and define the vector $\bv(\bx)=(v_1(\bx), \nv_1(\bx), \cdots, v_n(\bx), \nv_n(\bx))$ by
$$ 
\forall i\in\{1,\cdots,n\} \quad v_i(\bx)=x_i-\frac{1}{\sharp C_i(\bx)}\sum_{k\in C_i(\bx)} x_k \textrm{ and }
\nv_i(\bx)=\frac{ \frac{\rp+\rm}{2}-\s\nx_i }{(\rp-\rm)(\s \vee 1)} \, \rm .
$$
See figure 3.
\begin{figure}[h!]
  \centering
  \includegraphics[scale=0.9]{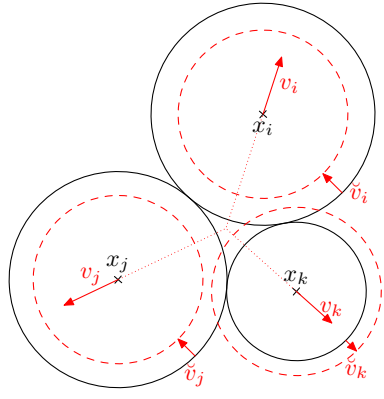}
  \caption{A cluster with 3 globules and the associated impulsion {\color{red}$\bv(\bx)$} constructed to push $\bx\in\partial\D$ back into the interior of the set of allowed globule configurations.}
\end{figure}

Since $|v_i(\bx)|\le(n-1)(2\rp)$ and $|\nv_i(\bx)|\le\frac{\rm}{2(\s \vee 1)}$, then \mbox{$|\bv(\bx)| \le 2 \rp n^{3/2}$.} Moreover~:
\begin{itemize}
\item
if $\bx\in\partial\D_{ij}$, then $C_i(\bx)=C_j(\bx)$ and
$$
\sqrt{2+2\s^2}~ \bv(\bx).\bn_{ij}(\bx)
= \s(\nx_i+\nx_j)-\s\rm\frac{\rp+\rm-\s\nx_i-\s\nx_j}{(\rp-\rm)(\s \vee 1)}
\ge 2\rm -\frac{\s\rm}{\s \vee 1} 
\ge \rm
$$
\item
if $\bx\in\partial\D_{i+}$, i.e. $\s\nx_i=\rp$, then $\bv(\bx).\bn_{i+}(\bx)=\frac{\rm}{2(\s \vee 1)}$
\item
if $\bx\in\partial\D_{i-}$, i.e. $\s\nx_i=\rm$, then $\bv(\bx).\bn_{i-}(\bx)=\frac{\rm}{2(\s \vee 1)}$.
\end{itemize}
So, with $\disp \beta_0= \frac{\rm}{4\rp(\s \vee 1) n^{3/2}}$, {\it (iv)} is satisfied.
\end{proof2}

\subsection{A full set of nice paths} \label{sec:nicepaths}
From now on, the techniques we use to study the globule model present some similarities with the methods developed for the model of hard balls treated in \cite{FR2}. So, in the rest of the paper, we will only detail the proofs which contain new technical difficulties. \\

We first bound from below the probability of globule paths which do not move too fast under the $(\En)$-dynamics.\\
For every $\eps>0$ and $\delta \in ]0,1]$, let $\tNP(\delta,\eps)$ denote the paths for which all globules  have a $\delta$-modulus of continuity
$\oscil$ smaller than $\eps$,  i.e.
$$
\tNP(\delta,\eps)
= \left\{  \bX \in \C([0,1],\A_g) :~ \forall i ,~ \oscil((X_i,\nX_i),\delta)\leq \eps  \right\} , 
$$
where the $\delta$-modulus of continuity of a globule path $(X,\nX)$ on $[0,1]$ is defined as
\begin{equation} \label{moduledecontinuite}
\qquad \qquad \qquad \quad \oscil((X,\nX),\delta)
:=\sup_{ \substack{0\le s,t\le 1\\ |t-s| \le \delta}} \sqrt{|X(t)-X(s)|^2 + (\nX(t)-\nX(s))^2} .
\end{equation} 
\begin{propal} \label{PropSpeedy}
There exists $c>0$ and $c_1 >0$ such that the following lower bound holds~:\\ 
$\forall \eps>0, \forall \delta \in ]0,1], \forall \ell \in \N^*$,                \\
$$
\inf_{\by \in \A_g} ~ Q^{\ell,\by} ( \tNP(\delta,\eps) ) ~~
\ge ~~1 - c_1 ~ \frac{\ell^3}{\delta} ~ \exp \left( - c \frac{\eps^2}{\delta} \right).
$$
\end{propal}

\begin{proof2} {\bf of Proposition \ref{PropSpeedy}}

By construction, the processes~
\begin{eqnarray*}
 W_i(t) &= & X^{\ell, \by,n}_i(t) - X^{\ell, \by,n}_i(0) \\
&& +\demi \int_0^t \nabla_x \psi^{\ell,\by}( X^{\ell, \by,n}_i(s),\nX^{\ell, \by,n}_i(s) ) ds 
   -\sum_{j=1}^n \int_0^t \frac{X^{\ell, \by,n}_i(s)-X^{\ell, \by,n}_j(s)}{\nX^{\ell, \by,n}_i(s)+\nX^{\ell, \by,n}_j(s)} dL_{ij}(s)\\
 \textrm{ and }\nW_i(t) &=& \frac{1}{\s}\big( \nX^{\ell, \by,n}_i(t) - \nX^{\ell, \by,n}_i(0)\big) \\
&& +  \frac{\s}{2} \int_0^t \nabla_{\nx} \psi^{\ell,\by}( X^{\ell, \by,n}_i(s),\nX^{\ell, \by,n}_i(s) ) ds 
   + \s \sum_{j=1}^n L_{ij}(t) + \frac{1}{\s}L_{i+}(t) - \frac{1}{\s}L_{i-}(t)
\end{eqnarray*}
are $3$-dimensional (resp. $1$-dimensional) Brownian motions starting from $0$.\\
When the initial distribution is $\nu^{\ell, \by}_n$ the law $Q_n^{\ell,\by}$ of $\bX^{\ell, \by,n}$ is reversible, and the backward processes
\begin{eqnarray*}
\widehat{W}_i(t)& =&  X^{\ell, \by,n}_i(1-t) - X^{\ell, \by,n}_i(1) \\
& &  + \demi \int_{1-t}^1 \nabla_x \psi^{\ell,\by}( X^{\ell, \by,n}_i(s),\nX^{\ell, \by,n}_i(s) )  ds 
     - \sum_{j=1}^n \int_{1-t}^1 \frac{X^{\ell, \by,n}_i(s)-X^{\ell, \by,n}_j(s)}{\nX^{\ell, \by,n}_i(s)+\nX^{\ell, \by,n}_j(s)} dL_{ij}(s)\\
  \textrm{ and } 
\widehat{\nW}_i(t)&  = & \frac{1}{\s}\big( \nX^{\ell, \by,n}_i(1-t) - \nX^{\ell, \by,n}_i(1)\big) 
                         + \frac{\s}{2} \int_{1-t}^1 \nabla_{\nx} \psi^{\ell,\by}( X^{\ell, \by,n}_i(s),\nX^{\ell, \by,n}_i(s) )ds \\
& & + \s \sum_{j=1}^n \big(L_{ij}(1)-L_{ij}(1-t)\big) - \frac{1}{\s}\big(L_{i+}(1)-L_{i+}(1-t)\big) + \frac{1}{\s}\big(L_{i-}(1)-L_{i-}(1-t)\big)
\end{eqnarray*}
are Brownian motions too.

As in \cite{LS}, the above equations provide the identities
$$
\left\{
\begin{array}{l}
\disp X^{\ell, \by,n}_i(t) - X^{\ell, \by,n}_i( 0) = \demi \Big(W_i(t)  + \widehat{W}_i(1-t)- \widehat{W}_i(1)\Big)  \\
\disp \nX^{\ell, \by,n}_i(t) - \nX^{\ell, \by,n}_i( 0) = \frac{\s}{2} \Big(\nW_i(t)  + \widehat{\nW}_i(1-t)- \widehat{\nW}_i(1)\Big) .
\end{array}
\right.
$$ 
Therefore, the control of the modulus of continuity of a globule path $(X^{\ell, \by,n}_i,\nX^{\ell, \by,n}_i)$ reduces to the estimate of the modulus of continuity of Brownian paths, as follows:
$$
\begin{array}{l}
Q_n^{\ell,\by} ( \tNP(\delta,\eps)^c  )                                      \\
 \disp
  \quad \le~ 2~ n ~ P \left( \oscil((W_1,\s \nW_1),\delta) > \eps \right) \, \nu_n^{\ell,\by}((\R^3\times \R)^n)\\
  \disp
  \quad \le~ 2~ n ~ P \left( \oscil((W_1,\s \nW_1),\delta) > \eps \right) \,  
\int_{\R^3\times [r_-,r_+]} \exp( - \psi^{\ell,\by}(x_1,\nx_1)) ~ ~dx_1 d\nx_1  ~ \nu_{n-1}^{\ell,\by}((\R^3\times \R)^{n-1})
\end{array}
$$
Now, we can use the following estimate obtained as corollary of Doob's inequality (for a proof in the one-dimensional case, see the Appendix of \cite{FR2})~:
\begin{lemma} \label{PropEstimeesBrownien}
Let us consider two independent Brownian motions  $ B \in \R^3$ and $\nB \in \R$.  
There exist two constants $c>0$ and $c_2>0$ (depending only on $\s$) such that for every $ \eps > 0 $ and every $ \delta \in ]0,1] $
$$
P ( \oscil((B,\s \nB),\delta) \ge \eps ) \quad \le \quad \frac{c_2}{\delta} \exp \left(- c\frac{\eps^2}{\delta} \right) 
$$
\end{lemma}
                                                                   \
This leads by summation  in $n$ to~:
$$
\begin{array}{l}
\disp 
Q^{\ell,\by}( \tNP(\delta,\eps)^c ) = \frac{e^{-|B(0,\ell)|}}{Z^{\ell,\by}} \sum_{n=1}^{+\infty} \frac{1}{n!} ~ Q_n^{\ell,\by} ( \tNP(\delta,\eps)^c  )  \\                                    \\
 \disp
  \quad \le~ 
 \frac{e^{-|B(0,\ell)|}}{Z^{\ell,\by}}  \Big(\sum_{n=1}^{+\infty} \frac{1}{(n-1)!} \nu_{n-1}^{\ell,\by}((\R^3\times \R)^{n-1}) \Big)~ \frac{2c_2}{\delta} \exp (- c\frac{\eps^2}{\delta} ) 
\int_{\R^3\times [r_-,r_+]} \exp( -\psi^{\ell,\by}(x_1,\nx_1)) ~ ~dx_1 d\nx_1  ~ \\
\disp
  \quad \le~ 
  \frac{2c_2 }{\delta} \exp(- c\frac{\eps^2}{\delta} ) \int_{\R^3\times [r_-,r_+]} \exp( -\psi^{\ell,\by}(x_1,\nx_1)) ~ ~dx_1 d\nx_1  .
\end{array}
$$
Recalling that $\psi^{\ell,\by}$ only vanishes into the ball $B(0,\ell)$, we get~:
\begin{eqnarray*}
&&\int_{\R^3\times [r_-,r_+]} \exp( - \psi^{\ell,\by}(x_1,\nx_1)) ~ ~dx_1 d\nx_1  \\
&\le & \int_{B(0,\ell)\times [r_-,r_+]} \exp( -\psi^{\ell,\by}(x_1,\nx_1)) ~ ~dx_1 d\nx_1  
 + \int \un_{\psi^{\ell,\by}>0} \exp( -\psi^{\ell,\by}(x_1,\nx_1)) ~ ~dx_1 d\nx_1  .
\end{eqnarray*}
The first term of the right side is smaller than $ \ell^3(r_+-r_-)|B(0,1)|$ and, thanks to (\ref{HypSurPsile}), the last term is uniformly bounded in $\ell$ and $\by$. 
This completes the proof of Proposition \ref{PropSpeedy}.
\end{proof2}
In order to control the convergence of the finite-dimensional systems, we have to estimate how many globules collide with a fixed globule $i$ during a short time interval. If the paths have a small oscillation, this set will be finite because globule $i$ can not reach  globules which are too far away.
But we also have to avoid the bump to propagate along a {\bf large chain} of neighboring globules.
We first define patterns called chains of globules, and then prove that they are rare enough, in the sense that their probability  decreases
exponentially fast as a function of the length of the chain. 
\begin{define}
Let $\eps>0$. The set of configurations containing an $\eps$-chain of M globules is defined by: 
$$ 
\Ch(\eps)~= \Big\{ \bx \in \A_g ,~\exists  i_1,\cdots,i_M \textrm{ distinct },~
                   |x_{i_2}-x_{i_1}|<\nx_{i_2}+\nx_{i_1}+\eps ,\cdots, |x_{i_M}-x_{i_{M-1}}|<\nx_{i_M}+\nx_{i_{M-1}}+\eps 
            \Big\}
$$
\end{define}
We now define a set of paths which are smooth in the sense that, at regular time intervals, there is no chain of globules: for $\delta \in1/\N^*,M\in\N^*, \eps>0$
$$
\ttNP(\delta,M,\eps)
:= \left\{ \bX \in \C([0,1],\A_g) :~\forall k \in \{0,\ldots,\frac{1}{\delta}-1\} ,~ \bX(\delta k) \not \in \Ch(\eps)  \right\}
$$
Note that this set decreases as a function of  $\eps$.\\
We now prove  a lower bound  for the $Q^{\ell,\by}$-Probability of $\ttNP(\delta,M,\eps)$.

\begin{propal} \label{PropBadChain}
For any $M \in \N^*$, there exists $c_3 >0$ such that, for any $\delta\in1/\N^* $ and  $0<\eps<1$~:
$$
\inf_{\by \in \A_g} ~ Q^{\ell,\by} \left( \ttNP(\delta,M,\eps) \right) ~\ge~ 1-\frac{c_3 }{\delta}~\ell^3~\eps^{M-1}.
$$
\end{propal}

\begin{proof2} {\bf of Proposition \ref{PropBadChain}}
\\
Let us estimate the $\nu_n^{\ell,\by}$- and the $\mu^{\ell,\by}$-Probability that a chain exists. For $n\ge M$:
\begin{eqnarray*}
\nu_n^{\ell,\by}( \Ch(\eps) ) &  \le &
\frac{n!}{(n-M)!} ~ \int_{(\R^3\times\R)^{n}} 
\prod_{i=2}^{M} \un_{\nx_i + \nx_{i-1} \le |x_i - x_{i-1}| \le \nx_i + \nx_{i-1}+\eps} ~ \prod_{i=1}^{M} \un_{r_- \le \nx_i \le r_+} \\
&&
\exp( - \sum_{i=1}^M \psi^{\ell,\by}(x_i,\nx_i)) \exp( - \sum_{i=M+1}^n \psi^{\ell,\by}(x_i,\nx_i))\\
&&
\prod_{i=M+1}^{n} \un_{r_- \le \nx_i \le r_+} 
\prod_{M+1\le i,j \le n} \un_{\nx_i + \nx_{j} \le |x_i - x_j| } 
 ~ dx_1 d\nx_1 \cdots dx_n d\nx_n  \\
&\le &
\frac{n!}{(n-M)!} \, \nu_{n-M}^{\ell,\by}((\R^3\times \R)^{n-M})~  \\
&&\int_{\R^3\times \R^M} |B(0,\nx_1 + \nx_{2}+\eps)\setminus B(0,\nx_1 + \nx_{2})|\cdots    
|B(0,\nx_{M-1} + \nx_{M}+\eps)\setminus B(0,\nx_{M-1} + \nx_{M})|   \\    
&& \prod_{i=1}^{M} \un_{r_- \le \nx_i \le r_+} ~ \exp( -\psi^{\ell,\by}(x_1,\nx_1)) ~dx_1  ~ d\nx_1 \cdots d\nx_M     \\
&\le &
   \frac{n!}{(n-M)!} \,  \nu_{n-M}^{\ell,\by}((\R^3\times \R)^{n-M})~ |B(0,2 r_+ +\eps)\setminus B(0,2r_+)|^{M-1}  \\
&&  (r_+-r_-)^{M-1} \int_{\R^3\times [r_-,r_+]} \exp( - \psi^{\ell,\by}(x_1,\nx_1)) ~ ~dx_1 d\nx_1  \\
& \le &
\frac{n!}{(n-M)!} \,  \nu_{n-M}^{\ell,\by}((\R^3\times \R)^{n-M})~  c_3 \ell^3 \eps^{M-1}
\end{eqnarray*}
for a certain constant $c_3>0$. Therefore
\begin{eqnarray*}
\mu^{\ell,\by}( \Ch(\eps) ) & \le&
\frac{e^{-|B(0,\ell)|}}{Z^{\ell,\by}} \sum_{n=M}^{+\infty} \frac{1}{(n-M)!} \nu_{n-M}^{\ell,\by}((\R^3\times \R)^{n-M})~  c_3 \ell^3 \eps^{M-1}\\
&\le & c_3  \ell^3 \eps^{M-1}
\end{eqnarray*}
and the stationarity of $Q^{\ell,\by}$ implies
$$
Q^{\ell,\by} \left( \ttNP(\delta,M,\eps)^c \right)
\le \sum_{k=0}^{\frac{1}{\delta}-1} Q^{\ell,\by}\left( \bX(\frac{k}{m}) \in \Ch(\eps) \right) \\ 
=  \frac{1}{\delta}  \, \mu^{\ell,\by}(\Ch(\eps))
\le \frac{c_3}{\delta} \ell^3 \eps^{M-1}.
$$
\end{proof2}
To prove the convergence of the approximations, we have to connect in a right way the different parameters $\delta$, $\eps$, $M$, $\ell$, in order to introduce a set of nice paths $\Omega_\by \subset \Omega$ on which the convergence holds. 
Since the Brownian motion has a.s. a $\delta$-modulus of continuity bounded by $\delta^\k$ for any $\k <1/2$, we choose $\k = 1/4$ and take $\eps$ proportional to $\delta^{1/4}$. \\
The maximal length $M$ of the chains is fixed, large enough so that $M>1+ \frac{4}{\k}=17$ (see (\ref{inequ:BadChain})).
Taking a unique scale parameter $m \in \N$ we thus choose
\begin{equation} \label{ChoixParametres}
\ell(m)= (1+3 r_+ )M 2^{4m}, \quad\quad \delta(m) = \frac{1}{2^{4m}}.
\end{equation}
It will be clear in (\ref{equ:Minorellm}) why this, with a suitable $\eps (m)=\frac{cst}{2^{m}}$, is a right choice. \\

We now define
\begin{eqnarray} \label{Omegax}
\Omega_\by &=& \liminf_{m \to +\infty}
               \left\{ \omega\in\Omega~: \bX^{\ell (m),\by}(\omega) \in 
                                         \tNP(\frac{1}{2^{4m}},\frac{1}{2^{m}}) \cap \ttNP(\frac{1}{2^{4m}},M,\frac{2^7}{2^{m}}) \right\}   .
\end{eqnarray}
We show in the next proposition that the set $\Omega_\by$ is of full measure with respect to any hard globule Poisson process.

\begin{propal} \label{PropNegligeabilite}
For any hard globule Poisson process $\mu\in\Pi_g$, one has
$$ \int_{\M} P(\Omega_\by) ~\mu(d\by) = 1 .$$
As a corollary,  for $\mu$ a.e. $\by, P(\Omega_\by)=1 $.
\end{propal}
\begin{proof2} We have to prove that $\disp \int_{\A_g} P(\Omega_\by^c) ~\mu(d\by) = 0$.

Thanks to Borel-Cantelli lemma, $\int_{\A_g} P(\Omega_\by^c)~\mu(d\by)$ vanishes as soon as the series 
\begin{eqnarray*}
& \sum_m & \int_{\A_g} P \Big(  \exists i: \oscil((X^{\ell(m),\by}_i,\nX^{\ell(m),\by}_i),\frac{1}{2^{4m}})> \frac{1}{2^m} \Big) \mu(d\by) \\
\textrm{ and } &\sum_m & \int_{\A_g}
P \Big( \exists k\leq 2^{4m} : \bX^{\ell(m),\by}(\frac{k}{2^{4m}}) \in \Ch(\frac{2^7}{2^{ m}}) \Big) \mu(d\by) 
\end{eqnarray*}
converges. 
Since for large $\ell, Q^{\ell,\by}$ and the law of $\bX^{\ell,\by}$ with initial distribution $\mu(\cdot|\by_{B(0,\ell)^c})$ are close, a similar argument as in \cite{FR2} Proof of Proposition 3.2 and the condition (\ref{HypSurPsile}) yield that these series converge as soon as
\begin{eqnarray}
& \sum_m & \int_{\A_g}
 Q^{\ell(m),\by} \Big(  \exists i: \oscil((X_i,\nX_i,\frac{1}{2^{4m}})> \frac{1}{2^m} \Big) \mu(d\by) < + \infty \label{ineq:oscpetites}\\
\textrm{ and } &\sum_m & \int_{\A_g}
Q^{\ell(m),\by} \Big( \exists k\leq 2^{4m} : \bX(\frac{k}{2^{4m}}) \in \Ch(\frac{2^7}{2^{m}}) \Big) \mu(d\by) < + \infty . \label{ineq:chainpetites}
\end{eqnarray}
Following Proposition \ref{PropSpeedy}
\begin{eqnarray*}
Q^{\ell(m),\by} \Big( \exists i: \oscil((X_i,\nX_i,\frac{1}{2^{4m}})>\frac{1}{2^{m}} \Big) 
&\le& c_1 \ell(m)^3 2^{4m} \exp\left( -c \, 2^{2m} \right)\\
&\le&  c_4 2^{16m} \exp \left( -c \, 2^{2m} \right)
\end{eqnarray*}
for a certain constant $c_4>0$ independent of $\by$. The above right side is the general term of a summable series in $m$. Therefore (\ref{ineq:oscpetites}) holds. \\
Following Proposition \ref{PropBadChain}
\begin{eqnarray} \label{inequ:BadChain}
Q^{\ell(m),\by} \Big( \exists k\leq 2^{4m} : \bX(\frac{k}{2^{4m}})  \in \Ch(\frac{2^7}{2^{m}}) \Big) 
&\le& c_3  ~ 2^{4m}  ~\ell(m)^3 ~ 2^{(7- m)(M-1)} \nonumber\\
&\le& c_5 2^{-m (M-17)}, 
\end{eqnarray}
for a certain constant $c_5>0$ independent of $\by$. We chose $M$ large enough to ensure the summability in $m$ of the right side, so (\ref{ineq:chainpetites}) holds and the proof is complete.
\end{proof2}

\subsection{The convergence} \label{sec:conv}

In this subsection, $\by\in\A_g$ is still fixed and we study the convergence of the approximating processes as $\ell \rightarrow +\infty $.  
\begin{propal} \label{PropCvgDesXl}
For every $\omega$ in $\Omega_\by$ and every $ i \in \N $, the sequence \\
$((X^{\ell(m),\by}_i, \nX^{\ell(m),\by}_i)(\omega,t), L_{i,j}^{\ell(m),\by}(\omega,t), L_{i+}^{\ell(m),\by}(\omega,t), L_{i-}^{\ell(m), \by}(\omega,t), j \in \N , t \in [0,1])_{m \in \N^*}$ is stationary as an element of $\C( [0,1], \R^3\times\R \times \R_+^{\N} \times \R_+^2 )$. 
The limit will be denoted by\\
$( (X^{\infty,\by}_i,\nX^{\infty,\by}_i)(\omega,t), L_{i,j}^{\infty,\by}(\omega,t), L_{i+}^{\infty,\by}(\omega,t), L_{i-}^{\infty,\by}(\omega,t), j \in \N , t \in [0,1] )$.
Therefore,
$$
\lim_{m \rightarrow  +\infty} \bX^{\ell(m),\by}(\omega,\cdot) = \bX^{\infty,\by}(\omega,\cdot)
$$ 
in $ \C([0,1], \A_g ) $.
\end{propal}
\begin{proof2}
The main idea is that if a fixed globule moves along a nice path, it will only  collide into a finite number of other globules. 
Thus dynamics $ ({\cal E}_g) $ reduces to an infinite number of SDE involving only a finite random number of particles up to time $1$.

Take $\omega \in \Omega_\by$ and $\rho>0$. Then, for $m$ large enough, $\bX^{\ell(m),\by}(\omega)$ and $\bX^{\ell(m+1),\by}(\omega)$ both belong to the same set of regular paths $\disp \tNP(\frac{1}{2^{4m}},\frac{2^3}{2^{m}}) \cap \ttNP(\frac{1}{2^{4m}},M,\frac{2^6}{2^{m}})$. 

For $\bX$ in this set and $k=0, \cdots, 2^{4m} -1$, we define the finite set of indices $J_{k,m}(\bX)$ as:
\begin{eqnarray*}
J_{k,m}(\bX) &:=& \Big\{ i\in\N, |X_i(\frac{k}{2^{4m}})| \le v_{k,m} \textrm{ or } \\
&& X_i(\frac{k}{2^{4m}}) \textrm{ belongs to some } \frac{2^6}{2^{m}}-\textrm{chain of globules which intersects }
B(0, v_{k,m}) \Big\},
\end{eqnarray*}
where $v_{k,m}:= \rho + (1+ 3r_+M)2^{4m} -3r_+M k$.\\
One can show (similarly as in Lemma 3.3 \cite{FR2}) that, for $m$ large enough:
\begin{eqnarray} \label{equ:ChaineInclusions}
\{i: |X_i(0)|\leq \rho \} \subset J_{2^{4m}-1,m}(\bX) \subset \cdots \subset J_{0,m}(\bX)
\end{eqnarray}
Moreover, a globule with index in $J_{k,m}(\bX)$ does not bump into globules outside this set~:
\begin{eqnarray} \label{equ:EcartCentres}
i \in J_{k,m}(\bX), j \not \in J_{k,m}(\bX) \Rightarrow 
\forall t \in [\frac{k}{2^{4m}},\frac{k+1}{2^{4m}} ] \quad |X_i(t) - X_j(t)| > \nX_i(t) + \nX_j(t) + \frac{2^5}{2^{4m}}
\end{eqnarray}
and it stays in a large ball around the origin~:
\begin{eqnarray} \label{equ:Minorellm}
i \in J_{k,m}(\bX) \Rightarrow &\forall t \in [\frac{k}{2^{4m}},\frac{k+1}{2^{4m}} ] &|X_i(t)|\leq v_{k-1,m} \le \ell(m)-2\rp .
\end{eqnarray}

Thus the penalization functions $\disp \psi^{\ell(m),\by}$ and $\disp \psi^{\ell(m+1),\by}$ vanish on globule $(X_i(t),\nX_i(t))$ if $i\in J_{k,m}(\bX)$ and $t\in [\frac{k}{2^{4m}},\frac{k+1}{2^{4m}} ]$.
Consequently, the paths $\bX^{\ell(m),\by}(\omega)$ and $\bX^{\ell(m+1),\by}(\omega)$ satisfy the following simplified version of equation ($\En$)~:
$$
\begin{array}{l}
 \forall k\in\{0, \cdots, 2^{4m} -1\}  \quad \forall i \in J_{k,m}(\bX) \quad 
 \forall t \in [\frac{k}{2^{4m}},\frac{k+1}{2^{4m}} ],                                                    \\ \disp
 X_i(t) 
 = X_i(\frac{k}{2^{4m}}) +W_i(\omega,t) -W_i(\omega,\frac{k}{2^{4m}})                                    
   + \sum_{j\in J_{k,m}(\bX)} \int_{\frac{k}{2^{4m}}}^t \frac{X_i(s)-X_j(s)}{\nX_i(s)+\nX_j(s)} dL_{ij}(s) \\ \disp
 \nX_i(t)
 = \nX_i(\frac{k}{2^{4m}}) +\s\nW_i(\omega,t) -\s\nW_i(\omega,\frac{k}{2^{4m}})                    \\ \disp\phantom{\nX_i(t) = }
   -\s^2 \sum_{j\in J_{k,m}(\bX)} \Big( L_{ij}(t)-L_{ij}(\frac{k}{2^{4m}})\Big)                       
   -\Big( L_{i+}(t)-L_{i+}(\frac{k}{2^{4m}})\Big) +\Big( L_{i-}(t)-L_{i-}(\frac{k}{2^{4m}}) \Big) .
\end{array}
$$
The initial configurations $\bX^{\ell(m),\by}(\omega,0)$ and $\bX^{\ell(m+1),\by}(\omega,0)$ are equal to the same configuration $\by$.
Hence the set of indices $J_{0,m}(\bX^{\ell(m),\by}(\omega))$ and $J_{0,m}(\bX^{\ell(m+1),\by}(\omega))$ are equal and 
$(X_i^{\ell(m),\by}(\omega),\nX_i^{\ell(m),\by}(\omega),i\in J_{0,m}(\bX^{\ell(m),\by}(\omega))$ satisfy the same equation as\\ 
\mbox{$(X_i^{\ell(m+1),\by}(\omega),\nX_i^{\ell(m+1),\by}(\omega),i\in J_{0,m}(\bX^{\ell(m+1),\by}(\omega))$} during the time interval $[0;\frac{1}{2^{4m}}]$.
The strong uniqueness in Proposition \ref{Prop:ExistApprox} implies the equality of the final values  $\bX_i^{\ell(m),\by}(\omega,\frac{1}{2^{4m}})$ and $\bX_i^{\ell(m+1),\by}(\omega,\frac{1}{2^{4m}})$ for the indices $i$ in the set $J_{0,m}(\bX^{\ell(m),\by}(\omega))$,
which contains both sets $J_{1,m}(\bX^{\ell(m),\by}(\omega))$ and $J_{1,m}(\bX^{\ell(m+1),\by}(\omega))$. 
Thus these two sets of indices are equal, which in turn implies that the paths $(X_i^{\ell(m),\by}(\omega),\nX_i^{\ell(m),\by}(\omega))$ and
$(X_i^{\ell(m+1),\by}(\omega),\nX_i^{\ell(m+1),\by}(\omega))$ coincide up to time $\frac{2}{2^{4m}}$ for indices $i\in J_{1,m}(\bX^{\ell(m),\by}(\omega))$. Using inclusions (\ref{equ:ChaineInclusions}) and the strong uniqueness again, we obtain the 
equality of both paths on $J_{2,m}(\bX^{\ell(m),\by}(\omega))=J_{2,m}(\bX^{\ell(m+1),\by}(\omega))$ up to time $\frac{3}{2^{4m}}$, and so on.

Strong uniqueness of the solution of ($\En$) holds for the path $\bX$ and the reflection term (linear combination of local times), but a priori not for each local time separately. However, as shown in the proof of corollary 3.6 in \cite{FradonBrownianGlobules}, the local times $L_{ij}$, $L_{i+}$, $L_{i}$ can be chosen in a unique way.  
With this choice, the same argument as above prove that local times $L_{i,j}^{\ell(m),\by}(\omega,t), L_{i+}^{\ell(m),\by}(\omega,t), L_{i-}^{\ell(m),\by}(\omega,t)$ and $L_{i,j}^{\ell(m+1),\by}(\omega,t), L_{i+}^{\ell(m+1),\by}(\omega,t), L_{i-}^{\ell(m+1),\by}(\omega,t)$ coincide for $i, j$ in $J_{0,m}(\bX^{\ell(m),\by}(\omega))$ and  $t \in [0;\frac{1}{2^{4m}}]$, and then again for $i, j$ in $J_{1,m}(\bX^{\ell(m),\by}(\omega))$ and $t \leq \frac{2}{2^{4m}}$, and so on.

In particular, if $|y_i| \le \rho$, then for $m$ large enough depending on $\rho$, $X_i^{\ell(m),\by}(\omega)=X_i^{\ell(m+1),\by}(\omega)$ and $\nX_i^{\ell(m),\by}(\omega)=\nX_i^{\ell(m+1),\by}(\omega)$ on the whole time interval $[0;1]$. The associated local times can be chosen in such a way that they also coincide, which implies the equality of the reflection terms up to time $1$. This completes the proof of the stationarity of the sequence of continuous functions $((X^{\ell(m),\by}_i, \nX^{\ell(m),\by}_i, L_{i,j}^{\ell(m),\by}, L_{i+}^{\ell(m),\by}, L_{i-}^{\ell(m), \by})(\omega,\cdot))_{m \in \N^*}$ and therefore its convergence to some path denoted by $(X^{\infty,\by}_i, \nX^{\infty,\by}_i, L_{i,j}^{\infty,\by}, L_{i+}^{\infty,\by}, L_{i-}^{\infty,\by}) (\omega,\cdot)$.\\
To check the convergence of $\bX^{\ell(m),\by}(\omega,\cdot)$ in $ \C([0,1], \A_g ) $, we remark that for each continuous function $f$ on $\R^3 \times \R $ with compact support,
$$
<\bX^{\ell(m),\by}(\omega,\cdot),f> = \sum_i f(X^{\ell(m),\by}_i(\omega,\cdot), \nX^{\ell(m),\by}_i(\omega,\cdot)) 
$$
where the sum is indeed finite due to the minimal distance $2r_-$ between any pair of points $X^{\ell(m),\by}_i(\omega,t)$ and $ X^{\ell(m),\by}_j(\omega,t)$ and the local boundedness of path oscillations. Therefore the stationary convergence of each term insures the stationary convergence of the sum.
\end{proof2}

\subsection{Properties of the limit process} \label{sec:limit}

To complete the proof of Theorem \ref{ThExistenceSolution} it suffices to show the following proposition.

\begin{propal} \label{PropUniciteLimite}
For every $\by \in \underline{{\A_g}} := \{ \bx \in \A_g: P(\Omega_\bx)=1 \}$ the family of processes \\
$( X^{\infty,\by}_i(t),\nX^{\infty,\by}_i(t), L_{i,j}^{\infty,\by}(t), L_{i+}^{\infty,\by}(t), L_{i-}^{\infty,\by}(t), i,j \in \N , t \in [0,1] )$ with initial configuration $\by$ solves uniquely the stochastic equation $(\E_g)$.
\end{propal}
\begin{proof2}
The equation satisfied by each 
$X^{\infty,\by}_i$ (resp. $\nX^{\infty,\by}_i$) includes by construction a finite sum of local time terms. It is straightforward to prove, using  a similar argumentation as in Proposition 4.1 of \cite{FR2}, that in fact this finite sum is already equal to the infinite sum present in $ (\E_g)$. 

The stationary convergence in Proposition \ref{PropCvgDesXl} implies that:
$$
\forall I \subset \N \text{ finite }, \exists m_0, \forall m\ge m_0 ~~~
(X^{\infty,\by}_i,\nX^{\infty,\by}_i)_{i\in I} \in \tNP(\frac{1}{2^{4m}},\frac{1}{2^{m}}) \cap \ttNP(\frac{1}{2^{4m}},M,\frac{2^7}{2^{m}})
$$
Consequently, using the strong uniqueness in Proposition \ref{Prop:ExistApprox} as in the proof of Proposition \ref{PropCvgDesXl}, uniqueness of the solution of $(\E_g)$ can be proved in the path space
$$
\Big\{ \bX \in \C([0,1], \A_g ):~ \forall I \subset \N \text{ finite }, \forall m_0, \exists m\ge m_0 ~~~
(X_i,\nX_i)_{i\in I} \in \tNP(\frac{1}{2^{4m}},\frac{1}{2^{m}}) \cap \ttNP(\frac{1}{2^{4m}},M,\frac{2^7}{2^{m}})
\Big\}
$$
\end{proof2}

Let us conclude with the proof of Proposition \ref{ThReversibiliteSolution}, that is with the reversibility of the solution of $ (\E_g)$ when the initial distribution is a hard globule Poisson process.
\begin{proof2} {\bf of Proposition \ref{ThReversibiliteSolution}}
Using similar estimates as in the proof of Proposition \ref{PropNegligeabilite}, the solution of $ (\E_g)$ starting with a hard globule Poisson process is the limit of processes whose distribution are close to  $Q^{\ell,\by}$, which is a time-reversible measure. More precisely, we have to prove that, if $\mu \in \Pi_g$,  for any $f_1,\ldots,f_k$ bounded continuous functions on $\A_g$ with compact support and for
$ t_1,\ldots,t_k \in [0,1]$
\begin{equation}
\label{DefReversibilite}
\int_{\A_g} \int_\Omega  \prod_{i=1}^k f_i(\bX^{\infty,\by}(\omega,t_i)) ~P(d\omega) ~\mu(d\by)
= \int_{\A_g} \int_\Omega  \prod_{i=1}^k f_i(\bX^{\infty,\by}(\omega,1-t_i)) ~P(d\omega) ~\mu(d\by)
\end{equation}
which is equivalent to 
$$
\lim_{m\to +\infty} 
\int_{\A_g} \int_\Omega  \Big( \prod_{i=1}^k f_i(\bX^{\ell(m),\by}(t_i)) 
- \prod_{i=1}^k f_i(\bX^{\ell(m),\by}(1-t_i)) \Big) ~dP ~\mu(d\by)
=0 .
$$
By computations similar to those done in \cite{FR2} to obtain inequality (17), we have\\
\begin{eqnarray*}
&&\Big|\int_{\A_g} \int_\Omega  \Big( \prod_{i=1}^k f_i(\bX^{\ell,\by}(t_i)) 
- \prod_{i=1}^k f_i(\bX^{\ell,\by}(1-t_i)) \Big) ~dP ~\mu(d\by)\Big|\\
&&\leq \Big|\int_{\A_g} \int_\Omega  \Big( \prod_{i=1}^k f_i(\bX(t_i)) 
- \prod_{i=1}^k f_i(\bX(1-t_i)) \Big) ~Q^{\ell,\by}(d \bX) ~\mu(d\by) \Big|\\
&& \quad + 2 ~ \prod_{i=1}^k \sup_{\bx \in \A_g} |f_i(\bx)| ~
            \int_{\A_g} \Big( 1 - \frac{Z^{B(0,\ell),\by}}{Z^{\ell,\by}} \Big) ~\mu(d\by)
\end{eqnarray*}
The first term of the right hand side is equal to 0.
The second term tends to zero as $\ell$ tends to infinity, thanks to assumption (\ref{HypSurPsile}).
\end{proof2}

\end{document}